\documentclass [a4paper,twoside,12pt]{article}

\usepackage[utf8]{inputenc}
\usepackage{dsfont}
\usepackage{amsmath,amsfonts,amssymb,color}
\usepackage{vmargin,graphicx,theorem}
\usepackage[english]{babel}
\usepackage{enumerate}
\usepackage{color}
\usepackage{pst-fill,pst-grad,pst-plot,pst-eucl,pstricks-add,pst-node}
\usepackage{todonotes}
\RequirePackage[colorlinks,linkcolor=blue,citecolor=blue,urlcolor=blue]{hyperref}

\setpapersize[portrait]{A4}
\setmarginsrb{1.5cm}{1cm}{1.5cm}{2.5cm}{1.5cm}{0cm}{0.5cm}{2cm}
\newcommand{\dR}{\ensuremath{\mathbf{R}}}

\newtheorem{ethm}{Theorem}[section]

\newtheorem{ecor}[ethm]{Corollary}

\newtheorem{eprop}[ethm]{Proposition}

\newtheorem{elem}[ethm]{Lemma}

\newtheorem{edefi}[ethm]{Definition}

\newtheorem{eex}[ethm]{Example}

\newcommand{\proofend}{~$\rhd$}
\newcommand{\proofbegin}{~$\lhd$}

\newenvironment{eproof}
               {\noindent {\emph{\textbf{Proof}}}\\\proofbegin~}
               {\proofend\\}

\renewcommand{\phi}{\varphi}
\renewcommand{\geq}{\geqslant}

\newcommand{\R}{\dR}

\newcommand{\e}{\varepsilon}

\newcommand{\grad}{\mathbf{grad}}
\newcommand{\supp}{\mathrm{supp}}

\newcommand{\beq}{\begin{equation}}\newcommand{\eeq}{\end{equation}}
\begin{document}

\title{Regularity properties of the Schrödinger cost}

\author{Gauthier Clerc\thanks{Institut Camille Jordan, Umr Cnrs 5208, Universit\'e Claude Bernard Lyon 1, 43 boulevard du 11 novembre 1918, F-69622 Villeurbanne cedex. Clerc@math.univ-lyon1.fr}}

\date{\today}
\maketitle

\abstract{The Schrödinger problem is an entropy minimisation problem on the space of probability measures. Its optimal value is a cost between two probability measures. In this article we investigate some regularity properties of this cost: continuity with respect to the marginals and time derivative of the cost along probability measures valued curves.}

\section{Introduction}

The Schrödinger problem was formulated by Schrödinger himself in the articles~\cite{Sch31, Sch32} in the thirties. The modern approach of this problem has been mainly developed in the two seminal papers~\cite{follmer1988} and~\cite{leonard2014}. The discovery in~\cite{mikami04} that the Monge-Kantorovitch problem is recovered as the short time limit of the Schrödinger problem has triggered intense research activities in the last decade. This interest is due to the fact that adding an entropic penalty in the Monge-Kantorovitch problem leads to major computationnal advantages using the Sinkhorn algorithm (see for instance~\cite{COTFNT}). The Schrödinger problem can also be a fruitful tool to prove some functional inequalities (see~\cite{clerc2020},~\cite{HWI}, etc.). \\  
The problem is, observing the empirical distribution of a cloud of Brownian particles at time $t=0$ and $t=1$, to find the distribution at time $0 < s < 1$ of the cloud. In modern language, this is an entropy minimization problem. The relative entropy of two measures is loosely defined by
$$
H(p|r):= \left\{\begin{array}{cc}
 \int  \log \left( \frac{dp}{dr}\right)dp \ \text{if} \ p \ll r ,\\
 + \infty\ \text{else}. 
\end{array} \right.
$$
We leave the precise definition of the relative entropy to the main body of the paper. Given two probability measures $\mu$ and $\nu$ on a Riemannian manifold $N$ equipped with a generator $L$ of reversible measure $m$, the Schrödinger cost is defined as
\begin{equation*}
\mathrm{Sch}(\mu,\nu):= \inf H(\gamma | R_{01}).
\end{equation*}
Here the infimum is taken over every probability measures on $N \times N$ with $\mu$ and $\nu$ as marginals and $R_{01}$ is the joint law of initial and final position of the unique diffusion measure with generator $L$ on $C\left( [0,1],N\right)$ starting from $m$. Independently proven in different papers (see~\cite{chengeorgiou2016, gentil-leonard2017, gentil-leonard2020, gigli-tamanini2020}), the Benamou-Brenier-Schrödinger formula state that 
\begin{equation*} \label{BBS1}
\mathrm{Sch}(\mu,\nu)= \frac{\mathcal{C}(\mu,\nu)}{4}+ \frac{H(\mu|m)+H(\nu|m)}{2}.
\end{equation*}
Here $\mathcal C(\mu,\nu)$ is the entropic cost given by
\begin{equation} \label{cout intro}
\mathcal{C}(\mu,\nu):= \inf \int_0^1 \left( \|v_s\|^2_{L^2(\mu_s)}+\left\|\nabla \log \frac{d \mu_s}{dm}\right\|_{L^2(\mu_s)}^2 \right) ds,
\end{equation}
and the infimum is taken over every $(\mu_s,v_s)_{0 \leq s \leq 1}$ such that $(\mu_s)_{0 \leq s \leq 1}$ is an absolutely continuous path with respect to the Wasserstein distance which connects $\mu$ to $\nu$ and satisfies in a weak sense for every $s \in [0,1]$
$$
\left\{\begin{array}{cc}
 v_s \in L^2(\mu_s),  \\
 \partial_s \mu_s=- \nabla \cdot (\mu_sv_s). 
\end{array} \right.
$$
In this paper we investigate regularity properties of the functions $(\mu,\nu) \mapsto \mathrm{Sch}(\mu,\nu)$ and $(\mu,\nu) \mapsto \mathcal C(\mu,\nu)$.{ To my knowledge, regularity properties of the Schrödinger cost as function over probability measures haven't been investigate yet, but the stability of optimizer has been investigate in~\cite{tamanini2017} and more recently in~\cite{ghosal2021stability} and~\cite{backhoff2019stability}.} We give an overview of the main contributions of this paper, leaving precise statements to other sections. 
\begin{itemize}
\item In Section~\ref{sec-cont} we investigate continuity properties of the cost function $\mathcal{C}$. In Theorem~\ref{continuityofthecost} we show that
$$
\underset{k \rightarrow \infty}{\lim} \mathcal{C}(\mu_k,\nu_k)= \mathcal{C}(\mu,\nu)
$$
if $W_2(\mu_k,\mu) \underset{k \rightarrow \infty}{\rightarrow} 0$ (resp. $\nu_k$ to $\nu$) with additional hypothesis about the entropy and the Fisher information along the sequences. 
\item In Section~\ref{sec-appli} we provide a few applications of the preceding continuity properties. The main result of this section is that using the continuity properties of $\mathrm{Sch}$ and $\mathcal{C}$ we are able to show that the Benamou-Brenier-Schrödinger formula~(\ref{BBS1}) is valid assuming that both measures have finite entropy, finite Fisher information and locally bounded densities. { Up to my knowledge this is a new result up to my knowledge, because no compactness assumptions are needed for the two marginals. }
\item In Section~\ref{sec-deriv} we investigate the question of the derivability of the functions $t \mapsto \mathrm{Sch}(\mu_t,\nu_t)$ and $t \mapsto \mathcal C (\mu_t,\nu_t)$, where $(\mu_t)_{t \geq 0}$ and $(\nu_t)_{t \geq 0}$ are some curves on the Wasserstein space. These results extend the existing ones for the Wasserstein distance, see~\cite[Theorem 8.4.7]{ambrosio-gigli2008} and~\cite[Theorem 23.9]{villani2009}. We prove that the derivative of the entropic cost is given for almost every $t$, by
$$
\frac{d}{dt} \mathcal C(\mu_t,\nu_t)= \langle \dot \mu_s^t|_{s=1}, \dot \nu_t \rangle_{L^2(\nu_t)} - \langle \dot \mu_s^t|_{s=0}, \dot \mu_t \rangle_{L^2(\mu_t)},
$$
where $(\mu_s^t)_{s \in [0,1]}$ is the minimizer of the problem~(\ref{cout intro}) from $\mu_t$ to $\nu_t$ and $\dot \mu_s^t$ is the velocity of the path $s \mapsto \mu_s^t$ defined in other section. Such minimizers are called entropic interpolations. Note that this is exactly the formula which holds for the Wasserstein distance, replacing the Wasserstein geodesics by the entropic interpolations. For technical reasons we prove this formula in the case where $N= \mathbb R^n$ and $L$ is the classical Laplacian operator.

\end{itemize}

\section{Setting of our work}

\subsection{Markov semigroups}

Let $(N, \mathfrak{g})$ be a smooth, connected and complete Riemannian manifold. We denote $dx$ the Riemannian measure and $\langle \cdot , \cdot \rangle$ the Riemannian metric (we omit $\mathfrak{g}$ for simplicity). Let $\nabla$ denote the gradient operator associated to $(N, \mathfrak{g})$ and $\nabla \cdot$ be the associated divergence in order to have for every smooth function $f$ and vector field $\zeta$
$$
\int \langle \nabla f(x) , \zeta(x) \rangle dx = -\int f(x) \nabla \cdot \zeta(x) dx.
$$
Hence the Laplace-Beltrami operator can be defined as $\Delta = \nabla \cdot \nabla$. We consider a differential generator $L:= \Delta - \langle \nabla , \nabla W \rangle $ for some smooth function $W:N \rightarrow \mathbb R$. We define the carré du Champ operator for every smooth functions $f$ and $g$ by
$$
\Gamma(f,g):= \frac{1}{2}\left(L(fg)-fLg-gLf \right).
$$
Under our current hypotheses we have $ \Gamma(f):= \Gamma(f,f)= |\nabla f|^2$, which is the length of $\nabla f$ with  to the Riemannian metric $\mathfrak{g}$. Let $Z:= \int e^{-W}dx$, then if $Z<\infty$ the reversible probability measure associated with $L$ is given by 
$$
d m:= \frac{e^{-W}}{Z}dx.
$$
If $Z= \infty$, the reversible measure associated with $L$ is $dm:= e^{-W}dx$ of infinite mass. Following the work of~\cite{bakry-emery1985} we define the iterated carré du champ operator given by
$$
\Gamma_2(f,g)= \frac{1}{2} \left( L \Gamma(f,g)-\Gamma(Lf,g)-\Gamma(f,Lg) \right), 
$$
for any smooth functions $f$ and $g$ and we denote $\Gamma_2(f):= \Gamma_2(f,f)$. We say that the operator $L$ verifies the $CD(\rho,n)$ curvature-dimension condition with $\rho \in \mathbb R$ and $n \in (0,\infty]$ if for every smooth function $f$
$$
\Gamma_2(f) \geq \rho \Gamma(f)+ \frac{1}{n}(Lf)^2.
$$
{  For instance, $\mathbb R^n$ endowed with the classical Laplacian operator verify the $CD(0,n)$ curvature-dimension condition. With the Ornstein-Uhlenbeck operator, $\mathbb R^n$ verify the $CD(1,\infty)$ curvature-dimension condition. More generally a Riemannian manifold of dimension $n \in \mathbb N$ and with a Ricci tensor bounded from below by $\rho \in \mathbb R$ endowed with his Laplace-Beltrami operator verify the $CD(\rho , n)$ curvature-dimension condition.} We assume that $L$ is the generator of a Markov semigroup $(P_t)_{t \geq 0}$, this is for example the case when a $CD(\rho,\infty)$ curvature-dimension condition holds for some $\rho \in \mathbb R$.  For every $f \in L^2(m)$ the family $(P_tf)_{t \geq 0}$ is defined as the unique solution of the Cauchy system
$$
\left\{ \begin{array}{cc}
 \partial_tu=Lu,  \\
 u(\cdot,0)=f( \cdot).
\end{array}
\right.
$$

Under the $CD(\rho,\infty)$ curvature-dimension condition this Markov semigroup admit a probability kernel $p_t(x,dy)$ with density $p_t(x,y)$, that is for every $t \geq 0$ and $f \in L^2(m)$
$$
\forall x \in N, \ P_tf(x)= \int f(y)p_t(x,dy)= \int f(y)p_t(x,y) dm(y), 
$$
for the existence of the kernel see~\cite[Theorem 7.7]{grigor}. We also define the dual semigroup $(P_t^*)_{t \geq 0}$ which acts on probability measures. Given a probability measure $\mu$ the family $(P_t^* \mu)_{t \geq 0}$ is given by the following equation
$$
\int f dP_t^{*} \mu = \int P_tf d\mu,
$$
 for every $t \geq 0$ and every test function $f$. When $\mu \ll m$, we have $\frac{dP_t^* \mu}{dm} =P_t \left( \frac{d \mu}{dm}e^W\right)$. The function $(t,x) \mapsto \frac{dP_t^* \mu}{dx}(x)$ is a solution of the following Fokker-Planck type equation
\begin{equation} \label{eq-51}
\partial_t \nu_t= L^* \nu_t := \Delta \nu_t+ \nabla \cdot \left(\nu_t \nabla W \right),
\end{equation}
with initial value $\frac{d \mu}{dx}$. Here $L^*$ is the dual operator of $L$ in $L^2(dx)$.

\subsection{Wasserstein space and absolutely continuous curves}

The set $\mathcal P_2(N)$ of probability measures on $N$ with finite second order moment can be endowed with the Wasserstein distance given for every $\mu,\nu \in \mathcal P_2(N)$ by,
$$
W_2^2(\mu,\nu):= \inf \sqrt{\int d^2(x,y)d \pi(x,y}),
$$
where the infimum is running over all $\pi \in \mathcal P(N \times N)$ with $\mu$ and $\nu$ as marginals { and $d$ is the Riemannian distance on (N, $\mathfrak{g}$)}. Recall that a path $(\mu_t)_{t \in [0,1]} \subset \mathcal P_2(N)$ is absolutely continuous with  to the Wasserstein distance $W_2$ if and only if
$$
|\dot \mu_t|:= \underset{s \rightarrow t}{\lim} \frac{W_2(\mu_s,\mu_t)}{|t-s|} \in L^1([0,1]).
$$
In this case, there exists a unique vector field $(V_t)_{t \in [0,1]}$ such that $V_t \in L^2(\mu_t)$ and $|\dot \mu_t|=\|V_t\|_{L^2(\mu_t)}$. Furthermore this vector field can be characterized as the solution of the continuity equation
$$
\partial_t \mu_t = - \nabla \cdot  (V_t \mu_t)
$$
with minimal norm in $L^2(\mu_t)$. We denote $\dot \mu_t=V_t$, and $(\dot \mu_t)_{t \in [0,1]}$ is called the velocity vector field of $(\mu_t)_{t \in [0,1]}$ or the velocity for short. Sometimes we also use the notation $\mathbf{dt} \mu_t= \dot \mu_t$. \\
In the famous paper~\cite{benamou-brenier2000} Benamou and Brenier showed that the Wasserstein distance admits a dynamical formulation
\begin{equation} \label{eq-50}
W_2^2(\mu,\nu)= \inf \int_0^1 \|\dot \mu_t\|^2_{L^2(\mu_t)}dt,
\end{equation}
where the infimum is running over all absolutely continuous paths which connect $\mu$ to $\nu$ in $\mathcal P_2(N)$. In his article~\cite{otto2001}, Felix Otto gave birth to a theory which allowed us to consider $(\mathcal P_2(N),W_2)$, heuristically at least, as an infinite dimensionnal Riemannian manifold. This theory was baptised \textquotedblleft Otto calculus\textquotedblright later by Cédric Villani. For every $\mu \in \mathcal P_2(N)$ the tangent space of $\mathcal P_2(N)$ at $\mu$ can be defined as
$$
T_{\mu} \mathcal P_2(N):= \overline{\left\{ \nabla \varphi: \ \varphi \in C_c^{\infty}(N)\right\}}^{L^2(\mu)},
$$
and the Riemannian metric is induced by the scalar product $\langle \cdot , \cdot \rangle_{L^2(\mu)}$, see for instance~\cite[Section 1.4]{gigli2012} or~\cite[Section 3.2]{gentil-leonard2017}. \\ 
As in the Riemannian case, the acceleration of a curve  can be defined as the covariant derivative of the veolcity field along the curve itself. If $(\mu_t)_{t \in [0,1]}$ is an absolutely continuous curve in $\mathcal P_2(N)$ and $(v_t)_{t \in [0,1]}$ is a vector field along $(\mu_t)_{t \in [0,1]}$, for every $t \in [0,1]$ we denote by $\mathbf{D_t}v_t$ the covariant derivative of $v_t$ along $(\mu_t)_{t \in [0,1]}$ defined in~\cite[Section 3.3]{gentil-leonard2017}. It turns out that in the case where the velocity field of $(\mu_t)_{t \in [0,1]}$ has the form $(\nabla \varphi_t)_{t \in [0,1]}$  then the acceleration of $(\mu_t)_{t \in [0,1]}$ is given by
$$
\forall t \in [0,1], \ \ddot \mu_t:=\mathbf{D_t} \dot \mu_t = \nabla \left( \frac{d}{dt} \varphi_t+ \frac{1}{2} \Gamma(\varphi_t)\right),
$$
see~\cite[Section 3.3]{gentil-leonard2017}. Covariant derivative and acceleration can be defined in more general framework, see~\cite[Section 5.1]{gigli2012}.

\subsection{Schrödinger problem} \label{sec-sch}

Here we introduce the Schrödinger problem by his modern definition, following the two seminal papers \cite{leonard2014} and \cite{follmer1988}. The first object of interest is the relative entropy of two measures. The relative entropy of a probability measure $p$ with  to a measure $r$ is loosely defined by
\begin{equation} \label{eq-41}
H(p|r):= \int \log \left( \frac{dp}{dr} \right)dp,
\end{equation}
if $p \ll r$ and $+ \infty$ elsewise. This definition is meaningful when $r$ is a probability measure but not necessarily when $r$ is unbounded. Assuming that $r$ is $\sigma$-finite, there exists a function $\mathcal{W}: M \rightarrow [1, \infty)$ such that $z_{\mathcal{W}}:= \int e^{-\mathcal{W}}dr < \infty$. Hence we can define a probability measure $r_{\mathcal{W}}:=z_{\mathcal{W}}^{-1}e^{-\mathcal{W}}r$ and for every measure $p$ such that $\int \mathcal{W} dp < \infty$
$$
H(p|r):=  H(p|r_{\mathcal{W}})- \int \mathcal{W} dp - \log(z_{\mathcal{W}}),
$$
where $H(p|r_{\mathcal{W}})$ is defined by the equation~(\ref{eq-41}). \\
For $\mu, \nu \in \mathcal P(N)$ we define the Schrödinger cost from $\mu$ to $\nu$ by
$$
\mathrm{Sch}(\mu,\nu):= \inf \left\{H(\gamma|R_{01}): \ \gamma \in \mathcal{P}(N \times N), \ \gamma_0=\mu, \ \gamma_1=\nu \right\},
$$
where $R_{01}$ is the joint law of the initial and final position of the Markov process associated with $L$ starting from $m$, which is given by
$$
dR_{01}(x,y)=p_1(x,y) dm(x)dm(y).
$$
To ensure the existence and unicity of minimizer, more hypothesis are needed. Namely we assume that there exists two non-negative mesurable functions $A,B:N \rightarrow \mathbb R$ such that
\begin{enumerate}[(i)] 
\item $p_1(x,y)  \geq e^{-A(x)-A(y)}$ uniformly in $x,y\in N$;
\item $\int e^{-B(x)-B(y)}p_1(x,y)\,m(dx) m(dy)< \infty$;
\item $\int (A+B)\, d\mu, \int (A+B)\, d\nu<\infty$;
\item $- \infty<H(\mu|m), H(\nu|m)<\infty$.
\end{enumerate}
We define the set
$$
\mathcal{P}^*_2(N):= \left\{ \mu \in \mathcal P_2(N): \ -\infty <  H(\mu|m) < \infty , \ \int \left(A+B \right) d\mu < \infty \right\}.
$$
If $\mu, \nu \in \mathcal P_2^*(N)$, it is proven that the Schrödinger cost $\mathrm{Sch}(\mu,\nu)$ is finite and admits a unique minimizer which takes the form
\begin{equation} \label{eq-57}
d\gamma= f \otimes g dR_{01},
\end{equation}
for two mesurable non-negative functions $f$ and $g$, see~\cite[Proposition 4.1.5]{tamanini2017}. Another fundamental result about the Schrödinger problem is an analogous formula to~(\ref{eq-50}) for the Schrödinger cost. 

\begin{ethm}[Benamou-Brenier-Schrödinger formula] \label{BBS-new} Let $\mu,\nu \in \mathcal P_2^*(N)$ be two probability measures compactly supported and with bounded densities with  to $m$. Then the following formula holds 
\begin{equation}
\label{eq-56}
 \mathrm{Sch}(\mu,\nu)=\frac{ \mathcal{C}(\mu,\nu)}{4}+\frac{\mathbf{\mathcal F}(\mu)+\mathbf{\mathcal F}(\nu)}{2},
\end{equation}
where $\mathcal{C}(\mu,\nu)$ is the entropic cost between $\mu$ and $\nu$ given by
$$
\mathcal{C}(\mu,\nu):= \inf \int_0^1 \left(\|\dot \mu_s\|_{L^2(\mu_s)}^2 + \mathcal \|\nabla \log (\mu_s)\|_{L^2(\mu_s)}^2\right) ds.
$$
Here the infimum is running over every absolutely continuous path $(\mu_s)_{s \in [0,1]}$ which connects $\mu$ to $\nu$ in $\mathcal P_2(N)$ and $\mathbf{\mathcal F}$ is defined as
$$
\mathbf{\mathcal F}(\mu):= H(\mu|m). 
$$
\end{ethm}
Different versions of this theorem have been obtained under various hypothesis, see~~\cite{chengeorgiou2016, gentil-leonard2017, gentil-leonard2020, gigli-tamanini2020}. \\
The functional $\mathbf{\mathcal F}: \mathcal P_2(N) \rightarrow [0,\infty]$ is central on this work. Its gradient can be identified by the equation $\frac{d}{dt} \mathbf{\mathcal F}(\mu_t)= \langle \mathbf{\grad}_{\mu_t} \mathbf{\mathcal F} , \dot \mu_t \rangle_{\mu_t}$ and is given for every $\mu \in \mathcal P_2(N)$ with smooth density against $m$ by
$$
\mathbf{\grad}_{\mu} \mathbf{\mathcal F} := \nabla \log \left(\frac{d \mu}{dm}\right). 
$$
Those definitions allowed us to see the Fokker-Planck type equation~(\ref{eq-51}) as the gradient flow equation of $\mathbf{\mathcal F}$. Indeed every solution $(\nu_t)_{t \geq 0}$ of this equation verify
$$
\dot \nu_t= - \nabla \left(\log \frac{\nu_t}{dx} +W \right) = - \nabla  \log \left(\frac{d \nu_t}{dm}\right)=- \mathbf{\grad}_{\nu_t} \mathbf{\mathcal F},
$$
see \cite[Section 3.2]{gentil-leonard2020}. With Otto calculus, we can also introduce the notions of Hessian and covariant derivative. A great fact is that the Hessian of $\mathbf{\mathcal F}$ can be expressed in term of $\Gamma_2$, indeed
$$
\forall \mu \in \mathcal P_2(N), \ \forall \ \nabla \varphi, \nabla \psi \in T_{\mu} \mathcal P_2(N), \ \mathbf{Hess}_{\mu} \mathbf{\mathcal F}(\nabla \varphi, \nabla \psi) = \int \Gamma_2(\nabla \varphi , \nabla \psi) d \mu,
$$
see~\cite[Section 3.3]{gentil-leonard2017}.
The quantity $\mathcal I(\mu):=\left\|\nabla \log \frac{d \mu}{dm}\right\|_{L^2(\mu)}^2$ which appears in the previous definition is central in this work, it is called the Fisher information. According to the Otto calculus formalism, the Fisher information admits the nice interpretation,
$$
\mathcal{I}(\mu):= \|\mathbf{\grad}_{\mu} \mathbf{\mathcal F} \|^2_{L^2(\mu)}. 
$$
Minimizers of the entropic cost $\mathcal{C}(\mu,\nu)$ are called entropic interpolations and take the form
$$
\mu_t=P_t f P_{1-t}g dm,
$$
where $f$ and $g$ are the two positive functions which appears in the equation~(\ref{eq-57}). Due to this particular structure, velocity and acceleration of entropic interpolations can be explicitly computed. It holds that for every $t \in [0,1]$
$$
\dot \mu_t= \nabla \left(\log  P_{1-t}g - \log P_t f\right).
$$
But the most important fact, is that entropic interpolations are solutions of the following Newton equation
\begin{equation*}
\ddot \mu_t= \nabla \frac{d}{dt} \log \mu_t + \nabla^2 \log \mu_t \dot \mu_t ,
\end{equation*}
which can be rewrite in the Otto calculus formalism as
\begin{equation} \label{eq-newton}
\ddot \mu_t = \mathbf{Hess}_{\mu_t} \mathbf{\mathcal F} \ \mathbf{\grad}_{\mu_t} \mathbf{\mathcal F}.
\end{equation}
This equation was first derived in~\cite[Theorem 1.2]{conforti2017}, see also~\cite[Sec 3.3, Propositon 3.5]{gentil-leonard2020}. 

\subsection{Flow maps} \label{sec-flowmaps}

In this subsection we follow~\cite[Sec 2.1]{gigli2012}. We need this result only in the euclidean framework, hence in this subsection we take $N= \mathbb R^n$ for simplicity. A crucial ingredient of the proof of the Theorem~\ref{lem-1} is, given a path $(\mu_t)_{t \in [0,1]}$, the existence of a family of maps $(T_{t \rightarrow s})_{t,s \in [0,1]}$ such that for every $s,t \in [0,1]$
$$
\frac{d}{ds}T_{t \rightarrow s}= \dot \mu_s \circ T_{t \rightarrow s}
$$
and
$$
T_{t \rightarrow s} \# \mu_t=\mu_s.
$$
These maps are called the flow maps associated with $(\mu_s)_{s \in [0,1]}$. The existence of such maps can be garanted by some regularity assumptions on the path. Before the statement we recall the definition of the Lipschitz constant of a vector field proposed by Gigli in~\cite{gigli2012}. 
\begin{edefi}[Lipschitz constant of a vector field] \label{constant}
For every smooth compactly supported vector field $\zeta$ on $\mathbb{R}^n$ we define
$$
L(\zeta):= \underset{x,y \in \mathbb R^n}{\sup} \frac{|\zeta(x)- \zeta(y)|}{|x-y|}
$$
Then for every $\mu \in \mathcal{P}_2(N)$ and every $v \in T_{\mu} \mathcal{P}_2(\mathbb R^n)$ we define
$$
\mathcal{L}(v):= \inf \underset{n \rightarrow \infty}{\underline{\lim}} L(\zeta_n),
$$
where the infimum is taken over sequences $(\zeta_n)_{n \in \mathbb N}$ of smooth compactly supported vector fields which converges to $v$ in $L^2(\mu)$ when $n \rightarrow \infty$.
\end{edefi}

Note that in the case where $v$ is smooth and compactly supported $\mathcal{L}(v)$ is the Lipschitz constant of $v$.

\begin{ethm}[Cauchy Lipschitz on manifolds,~\text{\cite[Theorem 2.6]{gigli2012}}]
Let $(\mu_t)_{t \in [0,1]} \subset \mathcal P_2(\mathbb R^n)$ be an absolutely continuous path such that
$$
\int_0^1 \mathcal{L}(\dot \mu_t)dt< \infty, \ \int_0^1 \| \dot \mu_t \|_{L^2(\mu_t)}dt < \infty.
$$
Then there exists a family of maps $(T_{t \rightarrow s})_{t,s \in [0,1]}$ such that
\begin{equation*}
\left\{ \begin{array}{cc}
T_{t \rightarrow s}: \supp(\mu_t) \rightarrow \supp (\mu_s),& \ \forall t,s \in [0,1],  \\
T_{t \rightarrow t} (x) =x,& \ \forall x \in \supp(\mu_t), \ t \in [0,1], \\
\frac{d}{dr}T_{{t \rightarrow r}_{|r=s}}= \dot \mu_s \circ T_{t \rightarrow s},& \ \forall t \in [0,1], \ a.e-s\in[0,1]. 
\end{array}
\right.
\end{equation*}
and the map $x \mapsto T_{t \rightarrow s}(x)$ is Lipschitz for every $s,t \in [0,1]$.
Furthermore for every $s,t,r \in [0,1]$ and $x \in \supp(\mu_t)$
$$
T_{r \rightarrow s} \circ T_{t \rightarrow r}(x)=T_{t \rightarrow s} (x),
$$
and
$$
T_{t \rightarrow s} \# \mu_t = \mu_s.
$$

\end{ethm}

\subsection{hypothesis about the heat kernel}

Here is a summary of all hypothesis needed in all the paper. 
\begin{enumerate}[(H1)]
\item The $CD(\rho,\infty)$ curvature-dimension condition holds for some $\rho \in \mathbb R$. 
\item hypothesis $(i)$ to $(iv)$ in Section~\ref{sec-sch}. 
\end{enumerate}
The first hypothesis~(H1) is needed to defined Markov semigroups as introduced in~\cite{bgl-book}. The second hypothesis~(H2) is needed to ensure existence and unicity of minimizers of the Schrödinger problem. For instance those hypothesis hold true when $N= \mathbb R^n$ is equipped with the classical Laplacian operator or the Ornstein-Ulhenbeck operator, or when $N$ is compact.

\section{Continuity of the entropic cost} \label{sec-cont}

Here we are interested in the continuity of the function $(\mu,\nu) \mapsto \mathcal{C}(\mu,\nu)$ where $\mathcal{C}(\mu,\nu)$ is defined as an infimum over all absolutely continuous paths connecting $\mu$ to $\nu$.

\begin{ethm}[Continuity of the entropic cost] \label{continuityofthecost} Let $\mu,\nu \in \mathcal{P}_2^*(N)$ and $(\mu_k)_{k \in \mathbb N},(\nu_k)_{k \in \mathbb N} \subset \mathcal P_2^*(N)$ be two sequences such that $\mu_k$ converges toward $\mu$ with  to the Wasserstein distance (resp.. $\nu_k$ toward $\nu$). We also assume that for every $k \in \mathbb N$ there exists an entropic interpolation from $\mu_k$ to $\nu_k$ (resp. from $\mu$ to $\nu$) and 
$$
\sup \left\{\mathcal{I}(\mu_k),\mathcal{I}(\nu_k); \, k \in \mathbb{N} \right\}< + \infty
$$
and 
$$
\sup \left\{\mathbf{\mathcal F}(\mu_k),\mathbf{\mathcal F}(\nu_k); \, k \in \mathbb{N} \right\}< + \infty.
$$
Then
$$
\mathcal{C}(\mu_k,\nu_k) \underset{k \rightarrow \infty}{\rightarrow} \mathcal{C}(\mu,\nu). 
$$
\end{ethm}

\begin{eproof}
To begin we will show that
$$
\underset{k \rightarrow \infty}{\overline{\lim}}\mathcal{C}(\mu_k,\nu_k) \leq \mathcal{C}(\mu,\nu). 
$$
To do so let us consider some particular path from $\mu_k$ to $\nu_k$. For every $k \in \mathbb{N}$, $\e \in (0,1/2)$ and $\delta \in \left(0, \e/2 \right)$, we define a path $\eta^{k,\e , \delta}$ from $\mu_k$ to $\mu$ given by 
$$
\eta_t^{k,\e,\delta}= \left\{ \begin{array}{cc}
 P_t^*(\mu_k), & \, t \in [0,\e/2-\delta),  \\
 P_{\e/2-\delta}^*(\gamma_t), & t \in [\e/2- \delta, \e/2+ \delta], \\
 P_{\e -t}^*(\mu), & t \in (\e/2+\delta, \e], 
\end{array}\right.
$$
where for all $t \in (\e/2- \delta, \e/2 + \delta)$ we define $\gamma_t= \alpha_{\frac{t-(\e/2 - \delta)}{2 \delta}}$ and $\left( \alpha_t\right)_{t \in [0,1]}$ is a Wasserstein constant speed geodesic from $\mu_k$ to $\mu$. We also define a path $(\tilde{\eta_t}^{k,\e,\delta})_{t \in [0,\e]}$ in the exact same way, but changing $\mu_k$ in $\nu$ and $\mu$ in $\nu_k$, that is a path from $\nu$ to $\nu_k$. \newline We denote by $\left(\mu_t\right)_{t \in [0,1]}$ the entropic interpolation from $\mu$ to $\nu$. Then for every $0 < \e < 1/2$, $k \in \mathbb{N}$ and $\delta \in (0, \e/2)$ we define a path $\left(\zeta_t^{k, \e, \delta}\right)_{t \in[0,1]}$ by
\begin{equation*}
\zeta_t^{k, \e, \delta}= \left\{
\begin{array}{cc}
 \eta^{k,\e,\delta}_t,& t \in [0,\e), \\
 \mu_{\frac{t-\e}{1-2 \e}},& t \in [\e,1-\e], \\
 \tilde{\eta}_{t-(1-\e)}^{k, \e, \delta},& t \in (1-\e,1].
\end{array}\right.
\end{equation*}
This is an absolutely continuous path which connects $\mu_k$ to $\nu_k$, hence, by the very definition of the cost $\mathcal C$ we have
$$
\mathcal{C}(\mu_k,\nu_k) \leq \int_0^1 \left\| \dot \zeta_t^{k,\e, \delta}\right\|^2_{L^2(\zeta_t^{k,\e, \delta})}+ \mathcal I(\zeta_t^{k, \e , \delta})dt.
$$
Due to the hypothesis (H1) we can apply the local logarithmic Sobolev inequalities stated in~\cite[Theorem 5.5.2]{bgl-book} and the $\rho$-convexity of the entropy (see~\cite[Corollary 17.19]{villani2009}) to find
$$
2 \int_{0}^{\e/2-\delta} \mathcal{I}(P_t^*\mu_k)dt+2 \int_{0}^{\e/2-\delta} \mathcal{I}(P_t^*\mu)dt \leq \frac{1-e^{- \rho\e}}{\rho} \left(\mathcal{I}(\mu)+ \mathcal{I}(\mu_k) \right),
$$
and
$$
\int_{\e/2-\delta}^{\e/2+\delta}\mathcal{I}(P_{\e/2-\delta}^* \gamma_t)dt \leq \frac{4 \delta \rho}{e^{\rho(\e-2\delta)}-1}\int_0^1 \mathbf{\mathcal F}(\alpha_t)dt \leq \frac{2 \delta \rho}{e^{\rho(\e-2\delta)}-1} \left( \mathbf{\mathcal F}(\mu) + \mathbf{\mathcal F} (\mu_k)- \frac{\rho}{2}W_2^2(\mu,\mu_k)\right).
$$
Here for $t \in (\e/2-\delta , \e/2+\delta)$ we denote $\mathbf{dt}P_{\e/2 - \delta} \gamma_t = \dot P_{\e/2 - \delta} \gamma_t$ the velocity field of the path $(P_{\e/2 - \delta} \gamma_s)_{s \in (\e/2 - \delta , \e/2+\delta)}$ at time $t$. We need to estimate
$$
\int_{\e/2 - \delta}^{\e/2 + \delta} \left\| \mathbf{dt} P_{\e/2-\delta}^*\gamma_t  \right\|_{L^2(P_{\e/2-\delta}^* \gamma_t)}dt.
$$
Using ~\cite[Theorem 8.3.1]{ambrosio-gigli2008}, for every $t \in (\e/2-\delta, \e/2+\delta)$ we have
\begin{equation*}
\left\| \mathbf{dt} P_{\e/2-\delta}^*\gamma_t  \right\|_{L^2(P_{\e/2-\delta}^*(\gamma_t))} = \underset{u \rightarrow t}{\lim} \frac{W_2(P_{\e/2- \delta }^*\gamma_t,P_{\e/2-\delta}^* \gamma_u)}{|t-u|}.
\end{equation*}
Finally, using the $CD(\rho,\infty)$ contraction property~\cite[Theorem 9.7.2]{bgl-book} we obtain \begin{equation*}\left\| \mathbf{dt} P_{\e/2-\delta}^*\gamma_t  \right\|_{L^2(P_{\e/2-\delta}^*(\gamma_t))} \leq e^{- \rho( \e/2 - \delta)} \frac{W_2(\mu_k,\mu)}{2 \delta}.\end{equation*}
We have shown 
\begin{multline*}
\int_0^{\e} \left\| \dot \zeta_t^{n,\e}\right\|^2_{L^2(\zeta_t^{n,\e})}+ \mathcal I\left(\zeta_t^{n,\e}\right)dt \leq \frac{1-e^{- \rho\e}}{\rho} \left(\mathcal{I}(\mu)+ \mathcal{I}(\mu_k) \right) \\ + \frac{2 \delta \rho}{e^{\rho(\e-2\delta)}-1} \left( \mathbf{\mathcal F}(\mu) + \mathbf{\mathcal F} (\mu_k)- \frac{\rho}{2}W_2^2(\mu,\mu_k)\right)  + e^{-\rho(\e - 2 \delta)} \frac{W_2^2(\mu_k,\mu)}{4\delta^2}.
\end{multline*}
A similar estimate hold for the integral from $1-\e$ to $1$ and we obtain
\begin{multline*}
\mathcal{C}(\mu_k,\nu_k)\leq \frac{1-e^{- \rho\e}}{\rho} \left(\mathcal{I}(\mu)+ \mathcal{I}(\mu_k) + \mathcal I(\nu) + \mathcal I(\nu_k) \right) + \\ \frac{2 \delta \rho}{e^{\rho(\e-2\delta)}-1} \left( \mathbf{\mathcal F}(\mu) + \mathbf{\mathcal F} (\mu_k) + \mathbf{\mathcal F}(\nu) + \mathbf{\mathcal F}(\nu_k)- \frac{\rho}{2}\left(W_2^2(\mu,\mu_k)+W_2^2(\nu,\nu_k)\right)\right) \\ + e^{-\rho(\e - 2 \delta)} \frac{W_2^2(\mu_k,\mu)+W_2^2(\nu_k,\nu)}{4\delta^2}+ \int_0^1 \frac{1}{1-2\e} \left\| \dot \mu_t \right\|_{L^2(\mu_t)}^2 + (1-2 \e)\mathcal{I}(\mu_t) dt.
\end{multline*}
Finally, letting in this order  $k$ tend to $\infty$, $\delta$ tend to $0$, and $\e$ tend to $0$ we obtain the desired inequality. 
\newline To obtain the $\liminf$ inequality, we consider the same path but swapping the role of $\mu_k$ and $\mu$ (resp. $\nu_k$ and $\nu$) and using the fact that $1-2 \e < \frac{1}{1-2\e}$, we obtain for every $k \in \mathbb N$, $\e \in (0,1/2)$ and $\delta \in (0,\e)$
\begin{multline*}
\mathcal{C}(\mu,\nu) \leq \frac{1-e^{-\rho \e}}{\rho}\left( \mathcal I(\mu) + \mathcal I(\mu_k)+ \mathcal I(\nu)+ \mathcal I(\nu_k)\right) \\ + \frac{2 \delta \rho}{e^{\rho(\e-2\delta)}-1} \left( \mathbf{\mathcal F}(\mu) + \mathbf{\mathcal F} (\mu_k) + \mathbf{\mathcal F}(\nu) + \mathbf{\mathcal F}(\nu_k)- \frac{\rho}{2}\left(W_2^2(\mu,\mu_k)+W_2^2(\nu,\nu_k)\right)\right) \\ + e^{-\rho(\e - 2 \delta)} \frac{W_2^2(\mu_k,\mu)+W_2^2(\nu_k,\nu)}{4\delta^2}+ \frac{1}{1-2 \e} \mathcal{C}(\mu_k,\nu_k).
\end{multline*}
Letting $k$ tends to $\infty$, $\delta$ tends to $0$ and $\e$ tend to zero we obtain
$$
\mathcal{C}(\mu,\nu) \leq \underset{k \rightarrow \infty}{\underline{\lim}}\mathcal{C}(\mu_k,\nu_k).
$$
\end{eproof}

\section{Extension of some properties to the non compactly supported case} \label{sec-appli}

\subsection{Benamou-Brenier-Schrödinger formula}

{  As mentionned before, the Benamou-Brenier-Schrödinger formula has been obtained under various hypothesis. Here we show that the result hold true in the case where both measures are not compactly supported but assuming that they have finite fisher information, finite entropy and locally bounded densities, using continuity properties of the cost proved before and existing results. Recall that, in the existing litterature, this formula is proved assuming that the two measures have bounded supports and densities, see~\cite[Theorem 4.3]{gigli-tamanini2020}}. 

\begin{eprop}[Benamou-Brenier-Schrödinger formula] \label{BBS-formula} Let $\mu,\nu \in \mathcal P_2^*(N)$ be two measures with locally bounded densities with respect to $m$ such that $\mathcal I(\mu), \mathcal I(\nu) < \infty$. Furthermore, assume that there exists an entropic interpolation from $\mu$ to $\nu$. Then $$
\mathrm{Sch}(\mu,\nu)= \frac{\mathcal{C}(\mu,\nu)}{4}+ \frac{\mathbf{\mathcal F}(\mu) + \mathbf{\mathcal F}(\nu)}{2}.
$$
\end{eprop}

Notice that, the hypothesis of existence of entropic interpolations is not so restrictive. Indeed if $N= \mathbb R^n$, entropic interpolations always exists for measures in $\mathcal{P}_2^*(N)$, see~\cite[Proposition 4.1]{leonard2014}.

\begin{eproof} Let $x \in N$, for every $n \in \mathbb N$, we define 
$$
\mu_n = \alpha_n \mathds{1}_{B(x,n)} \frac{d \mu}{d m} d m,
$$
where $\alpha_n$ is a constant renormalization. Analogously we can define a sequence $(\nu_n)_{n \in \mathbb N}$ which converges to $\nu$ when $n \rightarrow \infty$. As $\mu_n$ and $\nu_n$ are compactly supported, we can apply the Benamou-Brenier-Schrödinger formula, namely
\begin{equation} \label{eq-39}
\mathrm{Sch}(\mu_n,\nu_n)= \frac{\mathcal{C}(\mu_n,\nu_n)}{4}+ \frac{\mathbf{\mathcal F}(\mu_n) + \mathbf{\mathcal F}(\nu_n)}{2}.
\end{equation}
It can be easily shown that $W_2(\mu_n,\mu) \underset{n \rightarrow \infty}{\rightarrow} 0$, $\mathcal I(\mu_n) \underset{n \rightarrow \infty}{\rightarrow} \mathcal I(\mu)$, and $\mathbf{\mathcal F}(\mu_n) \underset{n \rightarrow \infty}{\rightarrow} \mathbf{\mathcal F}(\mu)$ (resp. $\nu_n$ and $\nu$). Hence by the theorem~\ref{continuityofthecost} the right-hand side of~(\ref{eq-39}) converges toward $\frac{\mathcal{C}(\mu,\nu)}{4}+ \frac{\mathbf{\mathcal F}(\mu) + \mathbf{\mathcal F}(\nu)}{2}$ when $n \rightarrow \infty$. \\

For the left hand-side, note that by the space restriction property of the Schrödinger cost~\cite[Proposition 4.2.2]{tamanini2017}, for every $n \in \mathbb N$ the optimal transport plan for the Schrödinger problem from $\mu_n$ to $\nu_n$ is given fro every probability set $A$ of $N \times N$ by
$$
\gamma_n(A):= \frac{\gamma(A \cap B(x,n)^2)}{\mu(B(x,n)) \nu(B(x,n))},
$$
where $\gamma$ is the optimal transport plan for the Schrödinger problem from $\mu$ to $\nu$. Hence $\mathrm{Sch}(\mu_n,\nu_n)=H(\gamma_n|R_{01}) \underset{n \rightarrow \infty}{\rightarrow} H(\gamma|R_{01})= \mathrm{Sch}(\mu,\nu)$, and the result is proved.
\end{eproof}

\subsection{Longtime properties of the entropic cost}

The entropic cost $\mathcal{C}(\mu,\nu)$ can be defined with more generality using a parameter $T>0$. For $\mu,\nu \in \mathcal P_2(N)$ and $T>0$ we define
$$
C_T(\mu,\nu):= \inf \int_0^T \|\dot \mu_t\|_{L^2(\mu_t)}+ \mathcal I (\mu_t) dt.
$$
In~\cite[Theorem 3.6]{clerclongtime2020} and~\cite[Theorem 1.4]{conforti2017}, estimates are provided for high values of $T$, but only in the case where both measures are compactly supported and smooth. Using the Proposition~\ref{continuityofthecost} we are able to extend these estimates to the non-compactly supported and non-smooth case. The following lemma will be very useful, it is proved in~\cite[Lemma 3.1]{HWI}.

\begin{elem}[Approximation by compactly supported measures] Let $\mu \in \mathcal P_2(N)$ be a probability measure such that $\mathbf{\mathcal F}(\mu) < \infty$ and $\mathcal I(\mu) < \infty$. Then there exists a sequence $(\mu_k)_{k \in \mathbb N} \subset \mathcal P_2(N)$ such that
\begin{enumerate}[(i)]
\item $\mathbf{\mathcal F}(\mu_k) \underset{k \rightarrow \infty}{ \rightarrow} \mathbf{\mathcal F}(\mu), \ \mathcal I(\mu_k) \underset{k \rightarrow \infty}{ \rightarrow} \mathcal I(\mu)$ and  $W_2(\mu_k,\mu) \underset{k \rightarrow \infty}{ \rightarrow} 0$.
\item $\frac{d \mu_k}{dm} \in C_c^{\infty}(N)$ for every $k \in \mathbb N$.
\end{enumerate}
\end{elem}

Using this lemma and the Theorem~\ref{continuityofthecost} we can easily extend the estimates provided in~~\cite[Theorem 1.4]{conforti2017} and~\cite[Theorem 3.6]{clerclongtime2020}. { Note that in~\cite{conforti2017} the author has already extended the estimate which holds under the $CD(\rho , \infty)$ curvature-dimension condition to the non-compact case, but we believe this is a pertinent example to illustrate the utility of Proposition~\ref{continuityofthecost}. The validity of the $CD(0,n)$ estimate for non-compactly supported measures is a new result up to my knowledge.}
 
\begin{ecor} [Talagrand type inequality for the entropic cost] Let $\mu,\nu \in \mathcal P_2(N)$ be two probability measures with finite entropy and Fisher information. Assume that there exists an entropic interpolation from $\mu$ to $\nu$. Then if the $CD(\rho,\infty)$ curvature-dimension condition holds for some $\rho >0$ 
$$
C_T(\mu,\nu) \leq 2 \underset{t \in (0,T)}{\inf} \left\{ \frac{1+e^{- 2 \rho t}}{1-e^{-2 \rho t}} \mathbf{\mathcal F}(\mu) + \frac{1+e^{- 2 \rho (T-t)}}{1-e^{-2 \rho (T-t)}} \mathbf{\mathcal F}(\nu) \right\}.
$$
If the $CD(0,n)$ curvature-dimension condition holds for some $n>0$ then 
$$
C_T(\mu,\nu) \leq C_1(\mu, \nu)+2n \log(T).
$$

\end{ecor}
These estimates are very useful, for instance they are fundamental to show the longtime convergence of entropic interpolations, see~\cite{clerclongtime2020}.

\section{Derivability of the Schrödinger cost} \label{sec-deriv}

In this section, we take $N= \mathbb R^n$ for some $n \in \mathbb N$ and $L= \Delta$ is the classical Laplacian operator. In this case the heat semigroup $(P_t)_{t \geq 0}$ is given by the following density
$$
\forall x,y \in \mathbb R^n, \ t >0, \ p_t(x,y)= \frac{1}{(4 \pi t)^{n/2}}e^{-\frac{|x-y|^2}{4t}},
$$
and the reversible measure $m$ is the Lebesgue measure. Notice that in this case, the funtions $A$ and $B$ wich appears in hypothesis $(i)$ to $(iv)$ in Section~\ref{sec-sch} can be chosen as
$$
\forall x \in \mathbb R^n, \ A(x)=B(x):=|x|^2.
$$
Hence in this case
$$
\mathcal P_2^*(\mathbb{R}^n)= \left\{ \mu \in \mathcal P_2(\mathbb R^n): \ - \infty < \mathcal{F}(\mu) < \infty \right\}.
$$

A natural question is the following: given a probability measure $\nu$ can we find a formula for the derivative of the function $t \mapsto \mathcal{C}(\mu_t,\nu)$ where $(\mu_t)_{t \in [0,1]}$ is a smooth curve in $\mathcal{P}_2(N)$? From a formal point view, we can easily find an answer. { Here we use the notation $\mathbf{dt} \mu_s^t$ (resp. $\mathbf{ds} \mu_s^t$)  for the velocity of a given path $(\mu_s^t)_{(s,t) \in [0,1] \times [0,1]}$ wrt to $t$ (resp. wrt $s$), to avoid confusion between the two variables.} For every $t \in [0,1]$, let $(\mu_s^t)_{s \in [0,1]}$ be the entropic interpolation from $\mu_t$ to $\nu$, then 
\begin{equation*}
\begin{split}
\frac{1}{2}\frac{d}{dt}\mathcal{C}(\mu_t,\nu)&=\frac{d}{dt} \int _0^1 \|\mathbf{ds}\mu_s^t\|^2_{L^2(\mu_s^t)}+\|\mathbf{\grad}_{\mu_s^t} \mathbf{\mathcal F}\|_{L^2(\mu_s^t)}^2ds \\
&=\int_0^1 \langle  \mathbf{D_t} \mathbf{ds} \mu_s^t, \mathbf{ds} \mu_s^t \rangle_{L^2(\mu_s^t)} + \mathbf{Hess}_{\mu_s^t} \mathbf{\mathcal F} ( \mathbf{ds} \mu_s^t, \mathbf{\grad}_{\mu_s^t} \mathbf{\mathcal F})ds \\
&= \int_0^1 \langle \mathbf{D_s} \mathbf{dt} \mu_s^t, \mathbf{ds} \mu_s^t \rangle_{L^2(\mu_s^t)} + \mathbf{Hess}_{\mu_s^t} \mathbf{\mathcal F} ( \mathbf{ds} \mu_s^t, \mathbf{\grad}_{\mu_s^t} \mathbf{\mathcal F })ds.
\end{split}
\end{equation*}
Here we have used~\cite[Lemma~20]{gentil-leonard2020} to invert the derivatives. Noticing that $$\langle \mathbf{D_s} \mathbf{dt} \mu_s^t, \mathbf{ds} \mu_s^t \rangle_{L^2(\mu_s^t)}= \frac{d}{ds} \langle \mathbf{dt} \mu_s^t, \mathbf{ds} \mu_s^t \rangle_{L^2(\mu_s^t)}- \langle \mathbf{dt} \mu_s^t , \mathbf{Ds} \mathbf{ds} \mu_s^t \rangle_{L^2(\mu_s^t) }$$ and using the Newton equation~(\ref{eq-newton}) we have
\begin{equation} \label{eq-1}
\frac{1}{2}\frac{d}{dt}\mathcal{C}(\mu_t,\nu)= \int_0^1  \frac{d}{ds} \langle \mathbf{dt} \mu_s^t, \mathbf{ds} \mu_s^t \rangle_{L^2(\mu_s^t)}ds= -\langle \dot \mu_t , \mathbf{ds} \mu_s^t|_{s=0} \rangle_{L^2(\mu_t)}.
\end{equation}
Unfortunately we do not see how to turn this proof into a rigorous one. \\

  From another point of view, we can try to derive the static formulation of the Schrödinger problem. Once again, we can easily guess a formula from a heuristic point of view. Indeed, let $(\mu_t)_{t \in [0,1]}$ be a smooth curve in $\mathcal P_2(N)$. For every $t \in [0,1]$ we denote by $\gamma_t=f^t \otimes g^t dR_{01}$ the optimal transport plan for the Schrödinger problem from $\mu_t$ to $\nu$. Then
\begin{equation*}
\begin{split}
\frac{d}{dt}\mathrm{Sch}(\mu_t,\nu)&=\frac{d}{dt}H(\gamma_t|R_{01}) \\
&= \langle \dot \gamma_t, \nabla \log \gamma_t \rangle_{L^2(\gamma_t)}. \\
\end{split}
\end{equation*}
Using the fact that $\gamma_t$ is a transport plan from $\mu_t$ to $\nu$ it can be easily shown that $\langle \dot \gamma_t, \nabla \log \gamma_t \rangle_{\gamma_t}=\langle \dot \mu_t, \nabla \log f^t  \rangle_{\mu_t}$. Hence we obtain
\begin{equation*}
\frac{d}{dt}\mathrm{Sch}(\mu_t,\nu)=\langle \dot \mu_t, \nabla \log f^t \rangle_{L^2(\mu_t)}.
\end{equation*}
Note that this is equivalent to the equation~(\ref{eq-1}) thanks to the Benamou-Brenier-Schrödinger formula. This proof is not rigorous because we don't have the regularity properties needed for $\gamma_t$. To prove our results, we follow the idea of Villani in~\cite[Theorem~23.9]{villani2009} where he computes the derivative of the Wasserstein distance along curves. Before the statement of our main theorem a technical lemma is needed. This lemma is an easy corollary of the proof of~\cite[Theorem 4.2.3]{tamanini2017}.

\begin{elem} \label{lemunif}
Let $(\mu_k)_{k \in \mathbb N},(\nu_k)_{k \in \mathbb N} \subset \mathcal P^*_2(\mathbb{R}^n)$ and $\mu,\nu \in \mathcal P^*_2(\mathbb{R}^n)$ such that $\mu_k$ converges toward $\mu$ with respect to the Wasserstein distance when $k \rightarrow \infty$ (resp. $\nu_k$ to $\nu$). For every $k \in \mathbb N$, we denote by $\gamma_k=f^k \otimes g^k dR_{01}$ the optimal transport plan for the Schrödinger problem from $\mu_k$ to $\nu_k$ and $\gamma = f \otimes g dR_{01}$ the optimal transport plan for the Schrödinger problem from $\mu$ to $\nu$. Assume that $(\frac{d \mu_k}{d m})_{k \in \mathbb N}$ and $(\frac{d \nu_k}{d m})_{k \in \mathbb N}$ are uniformly bounded in compact sets. Then for every compact set $K \subset N$, up to extraction $(f^k)_{k \in N}$ and $(g^k)_{k \in N}$ are uniformly bounded in $L^{\infty}(K,m)$. Furthermore
$$
f^k \underset{k \rightarrow \infty}{\overset{*}{\rightharpoonup}} f, \ g^k \underset{k \rightarrow \infty}{\overset{*}{\rightharpoonup}}  g,
$$
where the weak star convergence is understood in $L^{\infty}(K,m)$.
\end{elem}

In addition to this lemma, the following fact is central in our proof. Given two probability measures $p,r$ on $\mathbb R^n$ and a smooth enough function $\varphi: \mathbb R^n \rightarrow \mathbb R^n$, we have
$$
\frac{d \varphi \#p}{dm}=\frac{\frac{dp}{dm}}{|\det J_{\varphi}|} \circ \varphi^{-1},
$$
where $|\det J_{\varphi}|$ is the Jacobian determinant of $\varphi$. We often refer to this result as the Monge-Ampère equation or the Jacobian equation, see~\cite[Theorem 11.1]{villani2009} or \cite[Lemma 5.5.3]{ambrosio-gigli2008}. Using this equation, we obtain
\begin{equation} \label{eq-devent}
H(\varphi \#p|r)=H(p|r)- \int \log |\det J_{\varphi}|dp + \int \left( \log \frac{dr}{dm}-\log \frac{dr}{dm} \circ \varphi \right)dp,
\end{equation}
where $|\det J_{\varphi}|$ is the Jacobian determinant of $\varphi$. Given a curve $(\mu_t)_{t} \subset \mathcal P_2(N)$ and a measure $\nu \in \mathcal P_2(N)$, the idea of the following proof is to apply equation~(\ref{eq-devent}) with $r=R_{01}$, $p= \gamma_t$ is the optimal transport plan for the Schrödinger problem from $\mu_t$ to $\nu$ and $\varphi= T_{t \rightarrow s} \times \mathrm{Id}$ to bound from above $\mathrm{Sch}(\mu_s,\nu)$, and then let $s \rightarrow t$. 

\begin{ethm}[Derivative of the Schrödinger cost] \label{lem-1}
Let $(\nu_t)_{t \in (t_1,t_2)} \subset \mathcal P_2^* (\mathbb{R}^n)$ and $(\mu_t)_{t \in (t_1,t_2)} \subset \mathcal P_2^*(\mathbb{R}^n)$ be two absolutely continuous curves, for some $(t_1,t_2) \subset \mathbb R$. Furthermore assume that
\begin{enumerate}[(i)]
    \item For every $t \in (t_1,t_2)$ the measure $\mu_t$ has smooth bounded density against $m$. 
    \item There exists a constant $C>0$ such that for every $t \in (t_1,t_2)$ and $x \in \mathbb R^n$ we have $|\dot \mu_t (x)| \leq C(1+|x|)$. \label{h5}
    \item The sequence $\left( \frac{d \mu_t}{d m}\right)_{t \in (t_1,t_2)}$ and $\left( \frac{d \nu_t}{d m}\right)_{t \in (t_1,t_2)}$ are uniformly bounded in compact sets.
    \item The functions $t \mapsto \mathbf{\mathcal F}( \mu_t)$ and $t \mapsto \mathbf{\mathcal F} (\nu_t)$ are derivable and $\frac{d}{dt} \mathbf{\mathcal F}(\mu_t)= \langle \nabla \log \frac{d \mu_t}{dm} , \dot \mu_t \rangle_{L^2(\mu_t)}$ (resp. $\mathbf{\mathcal F}(\nu_t)$). 
    \item $\int_{t_1}^{t_2} \mathcal{L}( \dot \mu_t) dt < \infty$ and $\int_{t_1}^{t_2} \| \dot \mu_t \|_{L^2(\mu_t)} dt < \infty$, where $\mathcal{L}(\dot \mu_t)$ is defined at definition~\ref{constant}. 
    \item \label{h6} For every $t \in (t_1,t_2)$ the functions $f^t$ and $g^t$ are in $L^1(m)$, where $(f^t,g^t)$ is the unique solution in $L^{\infty}(m) \times L^{\infty}(m)$ of the Schrödinger system
    $$
    \left\{\begin{array}{cc} 
    \frac{d \mu_t}{dm}=f^t P_1 g^t,  \\
    \frac{d \nu_t}{dm}=g^t P_1 f^t. 
    \end{array} \right.
    $$
\end{enumerate}
Then the application $t \mapsto \mathrm{Sch}(\mu_t,\nu)$ is differentiable almost everywhere and we have for almost every $t \in (t_1,t_2)$
$$
\frac{d}{dt}\mathrm{Sch}(\mu_t,\nu_t)=\langle \dot \mu_t, \nabla \log f^t \rangle_{L^2(\mu_t)}+\langle \dot \nu_t , \nabla \log g^t \rangle_{L^2(\nu_t)}.
$$
Furthermore for almost every $t \in (t_1,t_2)$ this equality can be rewritten as
$$
\frac{d}{dt} \mathcal C (\mu_t,\nu_t)= \langle \dot \nu_t , \mathbf{ds} {\mu_s^t}_{|s=1} \rangle_{L^2(\nu_t)} - \langle \dot \mu_t , \mathbf{ds} {\mu_s^t}_{|s=0} \rangle_{L^2(\mu_t)},
$$
where $(\mu_s^t)_{s \in [0,1]}$ is the entropic interpolation from $\mu_t$ to $\nu_t$.
\end{ethm}

\begin{eproof}
To begin we want to show
$$
\frac{d }{dt}\mathrm{Sch}(\mu_t,\nu)=\langle \dot \mu_t , \nabla \log f^t \rangle_{L^2(\nu_t)},
$$
for every $\nu \in \mathcal P_2^*(\mathbb R^n)$ such that $\frac{d \nu}{dm} \in L^{\infty}(m)$. \\

For every $t \in [0,1]$, $\gamma_t$ denotes the optimal transport plan in the Schrödinger problem from $\mu_t$ to $\nu$. Let $t \in [0,1]$ be fixed. Then for every $s$ small enough by the very definition of the cost $\mathrm{Sch}(\mu_{t+s},\nu) \leq H((T_{t \rightarrow t+s} \times \mathrm{Id}) \# \gamma_t|R_{01})$ where $(T_{t_1 \rightarrow t_2})_{t_1,t_2 \in [0,1]}$ are the flow maps associated to $(\mu_s)_{s \in [0,1]}$ defined in the subsection~\ref{sec-flowmaps}. Applying the equation~(\ref{eq-devent}) with $r=R_{01}$, $p= \gamma_t$ and $\varphi = T_{t \rightarrow t+s} \times \mathrm{Id}$ we obtain
\begin{multline*}
H((T_{t \rightarrow t+s} \times \mathrm{Id}) \# \gamma_t |R_{01}) = H(\gamma_t |R_{01}) + \int \log  p_1 d \gamma_t- \int \log p_1  d (T_{t \rightarrow t+s} \times \mathrm{Id}) \# \gamma_t \\ -\int \log |\det J_{T_{t \rightarrow t+s}}(x)| d \mu_t(x).
\end{multline*}
As noticed in~\cite[Eq~(23.11)]{villani2009}, by the hypothesis~$(\ref{h5})$ there exists a constant $C$ such that for every $y \in \mathbb R^n$ ans $s_1,s_2 \in [0,1]$
\begin{equation} \label{eq-5}
\left\{ \begin{array}{cc}
 |T_{s_1 \rightarrow s_2}(x)| \leq C(1+|x|)  \\
 |x - T_{s_1 \rightarrow s_2}(x)| \leq C|s_1 - s_2|(1+|x|).
\end{array}\right.
\end{equation}
For every $x,y \in \mathbb R^n$, we have $\log p_1(T_{t \rightarrow t+s} x,y)=- \frac{|T_{t \rightarrow t+s}x-y|^2}{4}- \frac{n}{2} \log (4 \pi)$ and 
\begin{equation*}
 \begin{split}
 \left| \frac{1}{2}\frac{d}{ds}|T_{t \rightarrow t+s} (x)-y|^2 \right|&= \left| \langle \dot \mu_{t+s}\ \circ T_{t \rightarrow t+s}(x), T_{t \rightarrow t+s}(x)-y  \rangle \right| \\
 & \leq \frac{|\dot \mu_{t+s}\circ T_{t \rightarrow t+s}(x)|^2}{2} + \frac{|T_{t \rightarrow t+s}(x)-y|^2}{2} \\ 
 & \leq C(1+|x|^2+|y|^2) \in L^1(\gamma_t),
 \end{split}
\end{equation*}
for some constant $C>0$. Hence we can differentiate over the integral at time $s=0$ to find
\begin{equation} \label{eq-4}
\int \log  p_1 \ d (T_{t \rightarrow t+s} \times \mathrm{Id}) \# \gamma_t = \int \log  p_1  d\gamma_t   + s \int \langle \dot \mu_{t}(x), \nabla_x \log p_1(x,y)\rangle d \gamma_t(x,y)+o(s).
\end{equation}
Notice that thanks to the Monge-Ampère equation we have 
$$
\int \log |\det J_{T_{t \rightarrow t+s}} | d \mu_t= \mathbf{\mathcal F} (\mu_t) - \mathbf{\mathcal F} (\mu_{t+s}) = -s  \langle \nabla \log \frac{d\mu_t}{dm},  \dot \mu_t \rangle_{L^2(\mu_t)} + o(s).
$$
Combining this with the equation~(\ref{eq-4}), we have 
$$
\mathrm{Sch}(\mu_{t+s},\nu) \leq \mathrm{Sch}(\mu_t,\nu)-s \left( \int \langle \nabla_x \log p_1(x,y),\dot \mu_t(x) \rangle d \gamma_{t}(x,y) - \langle \nabla \log \frac{d\mu_t}{dm},  \dot \mu_t \rangle_{L^2(\mu_t)}  \right) + o(s). 
$$

Observe that using the hypothesis~$(\ref{h6})$ we have
\begin{equation*}
\begin{split}
\int \langle \nabla_x \log p_1(x,y),\dot \mu_t(x) \rangle d \gamma_{t}(x,y)&= \int \int \langle \nabla_x p_1(x,y),\dot \mu_t(x) \rangle f^t(x)g^t(y)dm(x)dm(y) \\
&= \int \langle \nabla_x \int p_1(x,y) g^t(y) dm(y) , \dot \mu_t (x) \rangle f^t(x) dm(x) \\
&= \int \langle \nabla P_1g^t(x), \dot \mu_t(x) \rangle f^t(x) dm(x) \\
&= \int \langle \nabla \log \frac{d \mu_t}{dm}(x), \dot \mu_t(x) \rangle d \mu_t(x) - \int \langle \nabla \log f^t(x), \dot \mu_t(x)\rangle d\mu_t(x).
\end{split}
\end{equation*}
Hence we obtain
$$
\overline{\underset{s \rightarrow 0}{\lim}}\frac{\mathrm{Sch}(\mu_{t+s},\nu)-\mathrm{Sch}(\mu_t,\nu)}{s} \leq \langle \dot \mu_t , \nabla \log f^t\rangle_{L^2(\mu_t)}.
$$

For the reverse inequality we use the same kind of estimates. By definition we have $\mathrm{Sch}(\mu_t,\nu) \leq H((T_{t+s \rightarrow t} \times Id) \# \gamma_{t+s}|R_{01})$. Applying equation~(\ref{eq-devent}) we have
\begin{multline*}
H((T_{t+s \rightarrow  t} \times Id) \# \gamma_{t+s}|R_{01})= H(\gamma_{t+s}|R_{01})- \int \log |\det J_{T_{t+s \rightarrow t}}| d \mu_{t+s} \\ - \int  \left( \frac{|T_{t+s \rightarrow t}x-y|^2-|x-y|^2}{4} \right) d \gamma_{t+s}.
\end{multline*}
As already noticed we have $\int \log  |\det J_{T_{t+s \rightarrow t}}| d \mu_{t+s}= \mathbf{\mathcal F}(\mu_{t+s}) - \mathbf{\mathcal F}(\mu_t)= s \langle \dot \mu_t , \nabla \log \mu_t\rangle_{L^2(\mu_t)}+o(s)$.
Now we have to deal with a more complicated term. We want to show that
$$
\int  \left( \frac{|T_{t+s \rightarrow t}x-y|^2-|x-y|^2}{4} \right) d \gamma_{t+s}(x,y)= s \int \langle \nabla_x \log  p_1 (x,y), \dot \mu_t(x) \rangle d \gamma_t(x,y) + o(s).
$$
Notice that using~(\ref{eq-5}) we have for every $s>0$
\begin{equation} \label{eq-6}
 \left||T_{t+s \rightarrow t}x-y|^2-|x-y|^2 \right| \leq  Cs (1+|x|^2+|y|^2)  
\end{equation}
for some $C>0$.

 For every $s \in \mathbb R$ small enough, we denote $v_s(x,y)= \frac{|T_{t+s \rightarrow t}x-y|^2-|x-y|^2}{s} $ and $v(x,y)=-2 \langle x-y , \dot \mu_{t}(x)\rangle$. Of course for every $x,y \in N$, we have $$v_s(x,y) \underset{s \rightarrow 0}{\rightarrow}v(x,y)$$ and by~(\ref{eq-6})
$$
| v_s(x,y) | \leq P(x,y) :=C(1+|x|^2+|y|^2).
$$
Let $\chi_R$ be the product function $\chi_R= \mathds{1}_{B(0,R)} \otimes \mathds{1}_{B(0,R)} $.
By the Lemma~\ref{lemunif}, for every $R>0$ there exists a sequence $(s_k^R)_{k \in \mathbb N}$ which tends to zero when $k$ tend to $\infty$ such that the sequences $(f^{t+s_k^R})$, $(g^{t+s_k^R})$ are uniformly bounded in $L^{\infty}(B(0,R),m)$ and 
$$
f^{t+s_k^R} \underset{k \rightarrow \infty}{\overset{*}{\rightharpoonup}}f^t,\ g^{t+s_k^R} \underset{k \rightarrow \infty}{\overset{*}{\rightharpoonup}}g^t,
$$
where the weak star convergence is understood in $L^{\infty}(K^R,m)$. Now for simplicity we denote $s_k^R=s_k$ and $K^R=B(0,R)$. \\

Note that
\begin{multline} \label{eq-7}
\int v_{s_k^R}(x,y) d \gamma_{t + s_k^R} (x,y) - \int v(x,y) d \gamma_t(x,y) = \int (1- \chi_R(x,y))v_{t+s_k^R}(x,y) d \gamma_{t+s_k^R}(x,y) \\ +  \int \chi_R(x,y) \left(v_{t+s_k^R}(x,y)f^{t+s_k^R}(x)g^{t+s_k^R}(y)-v(x,y)f^t(x)g^t(y) \right)dR_{01}(x,y) \\ + \int ( \chi_R(x,y)-1)v(x,y) d \gamma_t(x,y). \end{multline}
To obtain the desired estimate we are going to pass to the limsup in $k$, then let $R$ tend to $+ \infty$. The third term is independent of $k$ and by the dominated convergence theorem it is immediate that it tend to $0$ when $R \rightarrow \infty$. Things are trickier for the second term. Denote 
$$
\phi_k(x):= \int v_{s_k^R}(x,y)g^{t+s_k^R}(y) p_1(x,y) dm (y) 
$$
and 
$$
\phi(x):= \int v(x,y)g^{t}(y) p_1(x,y) dm (y). 
$$
Then
\begin{equation*}
\begin{split}
& \left| \int \chi_R(x,y)v_{s_k^R}(x,y) d \gamma_{t+s_k^R}(x,y) - \int \chi_R(x,y) v(x,y) d \gamma(x,y) \right| \\ &= \left| \int_{K} \phi_k f^{t+s_k^R}dm - \int_{K^R} \phi f^t dm \right| \\
& \leq \left| \int_{K^R} f^{t+s_k^R} (\phi^{k} - \phi) dm \right| + \left| \int_{K^R} \left( f^{t+ s_k^R} -f^t \right) \phi dm \right| \\
& \leq \underset{n \in \mathbb N}{\sup} \| f^{t+s_n}\|_{L^{\infty}(K^R,m)} \|\phi_k - \phi\|_{L^1(K^R,m)}+ \left| \int_{K^R} \left( f^{t+ s_k^R} -f^t \right) \phi dm \right|.
\end{split}
\end{equation*}
The second term tends to zero thanks to the weak star convergence of $f^{t+s_k^R}$ toward $f^t$ when $k \rightarrow \infty$. 
Furthermore the same kind of calculus gives for every $k \in \mathbb N$
\begin{multline*}
|\phi_k(x)-\phi(x)| \leq \underset{n \in \mathbb N}{\sup}\|g^{t+s_n}\|_{L^{\infty}(K,m)}  \|\left( v_k(x,\cdot) - v(x, \cdot)\right) p_1(x, \cdot)\|_{L^1(K,m)} \\ + \int v(x,y)(g(x)-g^{t+s_k^R}(x))dR_{01}(x,y).
\end{multline*}
Again the second term tends to zero thanks to the weak star convergence of $g^{t+s_k^R}$. Using the upper bound $$
\left|\left( v_k(x,\cdot) - v(x, \cdot)\right) p_1(x, \cdot) \right| \leq 2 P(x, \cdot)p_1(x, \cdot) \in L^1(K,m)
$$ we have by the dominated convergence theorem
$$
\|\left( v_k(x,\cdot) - v(x, \cdot)\right) p_1(x, \cdot)\|_{L^1(K,m)} \underset{k \rightarrow \infty}{\rightarrow 0}.
$$
Hence $\phi_k \underset{k \rightarrow \infty}{\rightarrow} \phi$ pointwise. Noticing that for every $x \in K^R$, we have $$|\phi_k(x)| \leq \underset{n \in \mathbb N}{\sup} \| g^{t+s_k^R}\|_{L^{\infty}(K^R,m)} \int P(x,y)p_1(x,y)dm(y) \in L^{1}(K^R,m).$$ 
By the dominated convergence theorem, $$ \|\phi_k-\phi\|_{L^1(K^R,m)} \underset{k \rightarrow \infty}{\rightarrow} 0.
$$
Thus the second term in~(\ref{eq-7}) tends to zero when $k \rightarrow + \infty$. For the first term term, notice that for $R \geq 1$
$$
\int (1- \chi_R) v_{t+s_k^R} d \gamma_{t+s_k^R} \leq 2C \int_{|x| \geq R}|x|^2 d \mu_{t+s_k^R}(x) + C \int_{|y| \geq R}|y|^2 d \nu(y),
$$
thus it converges to zero, see~\cite[Definition 6.8 and Theorem 6.9]{villani2009}.
Hence for every $R >0$ by letting $k$ tends to $+ \infty$ and $R$ tends to $+ \infty$ in
$$
\frac{Sch(\mu_t , \nu) - Sch(\mu_{t+s_k^R}, \nu)}{s_k^R} \leq - \frac{1}{s_k^R}\int \log |\det J_{T_{t+s_k^R \rightarrow t}}|d\mu_{t+s_k^R} + \int v_{s_k^R}(x,y) d \gamma_{t + s_k^R}(x,y)
$$ 
we obtain
$$
\underset{ s \rightarrow 0}{\overline{\lim}} \frac{Sch(\mu_t, \nu) - Sch(\mu_{t+s}, \nu)}{s} \leq \int \langle\nabla \log p_1(x,y),\dot \mu_t(x) \rangle d \gamma_{t}(x,y)+ \langle \dot \mu_t , \nabla \log \frac{d \mu_t}{dm} \rangle_{\mu_t}
$$
This is enough to conclude as in the previous case and obtain
$$
\underset{s \rightarrow 0}{\underline{\lim}}\frac{\mathrm{Sch}(\mu_{t+s},\nu)-\mathrm{Sch}(\mu_t,\nu)}{s} \geq \langle \dot \mu_t , \nabla \log f^t\rangle_{\mu_t}.
$$
This ends the case where $\nu_t= \nu$ is constant. Now we need to use a \textquotedblleft doubling of variables \textquotedblright technique. Let $s,s',t \in [0,1]$ and $\gamma_{s,t}$ (resp. $\gamma_{s',t}$) be the optimal transport plan for the Schrödinger problem from $\mu_s$ (resp. $\mu_s'$) to $\nu_t$. Then, using the same tricks as before we have
$$
H(\gamma_{s',t}|R_{01}) -H(\gamma_{s,t}|R_{01}) \leq \mathbf{\mathcal F}(\mu_{s'})- \mathbf{\mathcal F}(\mu_s)+ \frac{1}{4} \int \left( |x-y|^2 - |T_{s \rightarrow s'}x-y|^2 \right)d \gamma_{s,t}(x,y).
$$
Now using~(\ref{eq-6}), the fact that $s \mapsto \mathbf{\mathcal F} (\mu_s)$ is Lipschitz continuous and the fact that second order moment of of both curves are locally bounded, there exists a constant $C>0$ such that
$$
H(\gamma_{s',t}|R_{01}) -H(\gamma_{s,t}|R_{01}) \leq C|s-s'|.
$$
By symmetry we can take absolute values in this inequality and it follows that the function $(s,t) \mapsto \mathrm{Sch}(\mu_s,\nu_t)$ is locally absolutely continuous in $s$ uniformly in $t$ (also absolutely continuous in $t$ uniformly in $s$). Hence by~\cite[Lemma 23.28]{villani2009} the desired result follow.  
\end{eproof}

\begin{eex}[Contraction property along Gaussian curves] For $m, \sigma^2 >0$  we denote by $\mathcal{N}(m,\sigma^2)$ the probability measure on $\mathbb R$ given by
$$
d\mathcal{N}(m,\sigma^2)(x):= \frac{1}{\sqrt{2 \pi \sigma^2}} \exp\left(- \frac{(x-m)^2}{2 \sigma^2}\right)dx. 
$$
Considering the measures $\mu:= \mathcal{N}(m_0,1)$ and $\nu:= \mathcal{N}(m_1 , 1)$ it follows from~\cite[Sec A.2]{clerclongtime2020} that the curves $(P_t^* \mu)_{t \geq 0}$ and $(P_t^* \nu)_{t \geq 0}$ satisfies the hypothesis of Theorem~\ref{lem-1}. If we denote $(\mu_s^t)_{s \in [0,1]}$ the entropic interpolation from $P_t^* \mu$ to $P_t^* \nu$ applying Theorem 5.2 we obtain for almost every $t > 0$
\begin{equation} \label{eq-666}
  \begin{split}
  \frac{d}{dt} \mathcal C (P_t^* \mu , P_t^* \nu)&= - \langle \mathbf{\grad}_{P_t^* \nu} \mathcal Ent , \mathbf{ds} \mu_s^t|_{s=1}\rangle_{L^2(\nu_t)} + \langle \mathbf{\grad}_{P_t^* \mu} \mathcal Ent , \mathbf{ds} \mu_s^t|_{s=0}\rangle_{L^2(\mu_t)} \\
  &=-\frac{d}{ds} \mathcal Ent (\mu_s^t)|_{s=1} + \frac{d}{ds} \mathcal Ent(\mu_s^t)|_{s=0} \\
  &=- \int_0^1 \frac{d^2}{ds^2} \mathcal Ent (\mu_s^t) ds. \\
  \end{split}
\end{equation}
Using the Newton equation~(\ref{eq-newton}) and the $CD(0,1)$ curvature-dimension condition, for every $t >0 $ and $s \in [0,1]$ we have
\begin{equation*}
\begin{split}
\frac{d^2}{ds^2} \mathcal Ent(\mu_s^t)&= \mathbf{Hess}_{\mu_s^t} \mathcal Ent (\mathbf{ds} \mu_s^t, \mathbf ds \mu_s^t)+\mathbf{Hess}_{\mu_s^t}\mathcal Ent (\mathbf{\grad}_{\mu_s^t} \mathcal Ent, \mathbf{\grad}_{\mu_s^t} \mathcal Ent) \\
&\geq  \left( \left( \frac{d}{ds} \mathcal Ent(\mu_s^t)\right)^2 + |\mathbf{\grad}_{\mu_s^t}\mathcal Ent(\mu_s^t)|_{L^2(\mu_s^t)}^4\right) .
\end{split}
\end{equation*}
Then by the Jensen inequality and neglecting the second term we obtain for every $t \geq 0$,
$$
\int_0^1 \frac{d^2}{ds^2} \mathcal Ent(\mu_s^t) ds \geq  \left(\mathcal Ent ( P_t^* \mu) - \mathcal Ent (P_t^* \nu ) \right)^2.
$$
Using this and integrating the equality~(\ref{eq-666}) we obtain for every $t \geq 0$,
$$
\mathcal{C} (P_t^* \mu, P_t^* \nu) \leq \mathcal{C}(\mu,\nu) - \int_0^t \left(\mathcal Ent ( P_u^* \mu) - \mathcal Ent( P_u^* \nu ) \right)^2 du.
$$
Hence we recover the $(0,1)$-contraction property of the entropic cost proved in~\cite[Theorem 37]{gentil-leonard2020}. The result could, of course, be proven in $\R^n$ for $n\geq1$. 

\end{eex}

\paragraph{Acknowledgements} This work  was supported by the French ANR-17-CE40-0030 EFI project. { I want to thanks warmly the anonymous referees for their work. }

\bibliographystyle{alpha}
{\footnotesize{\bibliography{reg}}}
\end{document}